%&bigtex
\magnification=\magstep1
\input amstex
\voffset=-3pc
\documentstyle{amsppt}
\NoBlackBoxes
\def\Bp{\Bbb P}
\def\Bn{\Bbb N}
\def\Bc{\Bbb C}
\def\Br{\Bbb R}
\def\Bz{\Bbb Z}

\def\Co{\Cal O}
\def\Ce{\Cal E}
\def\supp{\text{supp}}

\def\const{\text{const}}
\def\Ker{\text{Ker}}
\def\SO{\text{SO}}
\topmatter
\title Dolbeault cohomology of a loop space\endtitle
\author
L\'aszl\'o Lempert and Ning Zhang\footnote"*"{Research of both authors 
was partially supported by an NSF grant.\hfill\hfill}
\endauthor
\address
Department of Mathematics, Purdue University, West Lafayette, IN 47907
USA
\endaddress
\address Department of Mathematics, University of California, 
Riverside,
CA 92521 USA
\endaddress
\abstract
The loop space $L\Bbb P_1$ of the Riemann sphere is an infinite 
dimensional
complex manifold consisting of maps (loops) $S^1\to \Bbb P_1$ in some 
fixed
$C^k$ or Sobolev $W^{k,p}$ space.
In this paper we compute the Dolbeault cohomology groups 
$H^{0,1}(L\Bp_1)$.
\endabstract
\endtopmatter
\document
\leftheadtext{L\'aszl\'o Lempert and Ning Zhang}
\rightheadtext{Dolbeault Cohomology of a Loop Space}
\hsize=6.5truein
\vsize=9truein
\parskip=6pt
\heading{0.\ Introduction}\endheading

Loop spaces $LM$ of compact complex manifolds $M$ promise to have rich 
analytic
cohomology theories, and it is expected that sheaf and Dolbeault 
cohomology
groups of $LM$ will shed new light on the complex geometry and analysis 
of $M$
itself.  This idea first occurs in [W], in the context of the infinite 
dimensional
Dirac operator, and then in [HBJ] that touches upon Dolbeault groups of 
loop spaces;  but in all this both works stay heuristic.
Our goal here is rigorously to compute the $H^{0,1}$ Dolbeault group of 
the first interesting loop space, that of the Riemann sphere $\Bbb 
P_1$.  The consideration of $H^{0,1}(L\Bbb P_1)$ was directly motivated by 
[MZ], that
among other things features a curious line bundle on $L\Bbb P_1$. More
recently, the second named author in [Z] classified all holomorphic
line bundles on $L\Bbb P_1$ that are invariant under a certain group
of holomorphic automorphisms of $L\Bbb P_1$---a problem closely related
to describing (a certain subspace of) $H^{0,1}(L\Bbb P_1)$.
One noteworthy fact that emerges from the present
research is that analytic cohomology of loop spaces, unlike topological
cohomology (cf. [P, Theorem 13.14]), is rather sensitive to the 
regularity
of loops admitted in the space.
Another concerns local functionals, a notion from theoretical physics.
Roughly, if $M$ is a manifold, a local functional on a space of loops 
$x\colon S^1\to M$ is one of form
$$
f(x)=\int_{S^1}\Phi(t,x(t),\dot x(t), \ddot x(t),\ldots)dt,
$$
where $\Phi$ is a function on $S^1\times$ \!an appropriate jet bundle 
of $M$.
It turns out that all cohomology classes in $H^{0,1}(L\Bp_1)$ are given 
by local functionals.
Nonlocal cohomology classes exist only perturbatively, i.e., in a 
neighborhood of constant loops in $L\Bp_1$; but none of them extends to the 
whole of $L\Bp_1$.

We fix a smoothness class $C^k$, $k=1,2,\dots,\infty$, or Sobolev 
$W^{k,p}$,
$k=1,2,\dots,1\le p<\infty$. If $M$ is a finite dimensional complex
manifold, consider the space $LM=L_k M$,
resp. $L_{k,p}M$ of maps $S^1=\Bbb R/\Bbb Z\to M$
of the given regularity.
These spaces are complex manifolds modeled on a Banach space, except
for $L_\infty M$, which is modeled on a Fr\'echet space. We
shall focus on the loop space(s) $L\Bbb P_1$. As on any
complex manifold, one can consider the space $C_{r,q}^\infty(L\Bbb 
P_1)$ of
smooth $(r,q)$ forms, the operators 
$\bar\partial:C_{r,q}^\infty(L\Bp_1)\to
C_{r,q+1}^\infty(L\Bp_1)$, and the associated Dolbeault groups 
$H^{r,q}(L\Bp_1)$;
for all this, see e.g. [L1,2].
On the other hand, let $\frak F$ be the space of holomorphic functions 
$F\colon\Bc\times L\Bc\to\Bc$ that have the following properties:
\roster
\item"(1)"$F(\zeta/\lambda,\lambda^2 y)=O(\lambda^2)$, as 
$\Bc\ni\lambda\to 0$;
\item"(2)"$F(\zeta,x+y)=F(\zeta,x)+F(\zeta,y)$, if $\supp\, x\cap\supp 
\, y=\emptyset$;
\item"(3)"$F(\zeta,y+\text{const})=F(\zeta,y)$.
\endroster
As we shall see, the additivity property (2) implies $F(\zeta,y)$ is 
local in $y$.

\proclaim{Theorem 0.1}$H^{0,1}(L\Bp_1)\approx\Bc\oplus\frak F$.
\endproclaim

In the case of $L_\infty\Bp_1$, examples of $F\in\frak F$ are
$$
F(\zeta,y)=\zeta^\nu\biggl<\Phi,\ \prod^m_{j=0} 
y^{(d_j)}\biggr>,\tag0.1
$$
where $\Phi$ is a distribution on $S^1, y^{(d)}$ denotes $d$'th 
derivative, each $d_j\geq d_0=1$, and $0\leq\nu\leq 2m$.
A general function in $\frak F$ can be approximated by linear 
combinations of functions
of form (0.1), see Theorem 1.5.

On any, possibly
infinite dimensional complex manifold $X$ the space
$C_{r,q}^\infty(X)$ can be given the compact--$C^\infty$ topology as 
follows.
First, the compact--open topology on $C_{0,0}^\infty(X)=C^\infty(X)$ is
generated by $C^0$ seminorms $\Vert f\Vert_K=\sup_K|f|$ for all 
$K\subset X$
compact.  The family of $C^\nu$ seminorms is defined inductively:  each
$C^{\nu-1}$ seminorm $\Vert\ \ \Vert$ on $C^\infty(TX)$ induces a 
$C^\nu$
seminorm $\Vert f\Vert'=\Vert df\Vert$ on $C^\infty(X)$.  The 
collection of
all $C^\nu$ seminorms, $\nu=0,1,\dots$, defines the compact--$C^\infty$ 
topology on $C^\infty(X)$.  The compact--$C^\infty$ topology on a 
general
$C_{r,q}^\infty(X)$ is induced by the embedding 
$C_{r,q}^\infty(X)\subset C^\infty(\overset{r+q}\to\oplus\ TX)$. With this topology 
$C_{r,q}^\infty(X)$ is
a separated locally convex vector space, complete if $X$ is first 
countable.
The quotient space $H^{r,q}(X)$ inherits a locally convex topology, not
necessarily separated.
We note that on the subspace $\Co(X)\subset C^\infty(X)$ of holomorphic 
functions the compact--$C^\infty$ topology restricts to the 
compact--open topology.
The isomorphism in Theorem 0.1 is topological; it is also equivariant 
with
respect to the obvious actions of the group of $C^k$ diffeomorphisms of 
$S^1$.

There is another group, the group $G\approx\text{PSL}(2,\Bc)$ of 
holomorphic automorphisms of $\Bp_1$, whose holomorphic action on $L\Bp_1$ 
(by
post--composition) and on $H^{0,1}(L\Bp_1)$ will be of greater concern 
to us.
Theorems 0.2, 0.3, 0.4 below will describe the structure of 
$H^{0,1}(L\Bp_1)$ as a $G$--module.
Recall that any irreducible (always holomorphic) $G$--module is 
isomorphic, for some $n=0,1,\ldots$, to the space $\frak K_n$ of holomorphic 
differentials $\psi(\zeta)(d\zeta)^{-n}$ of order $-n$ on $\Bp_1$; here 
$\psi$ is a polynomial, deg $\psi\leq 2n$, and $G$ acts by pullback.
(For this, see [BD, pp.~84-86], and note that the subgroup 
$\approx\SO(3)$
formed by $g\in G$ that preserve the Fubini--Study metric is a 
maximally real submanifold; hence the holomorphic representation theory of $G$ 
agrees with the representation theory of $\text{SO}(3)$.)
The $n$'th isotypical subspace of a $G$--module $V$ is the sum of all 
irreducible submodules isomorphic to $\frak K_n$.
In particular, the $0$'th isotypical subspace is the space $V^G$ of 
fixed vectors.

\proclaim{Theorem 0.2}If $n\geq 1$, the $n$'th isotypical subspace of 
$H^{0,1}(L_\infty\Bp_1)$ is isomorphic to the space $\frak F^n$ spanned 
by functions of form (0.1), with $m=n$.
\endproclaim

The isomorphism above is that of locally convex spaces, as $\frak F$ or
$\frak F^n$ have not been endowed with an action of $G$ yet.
But in Section 2 they will be, and we shall see that the isomorphism in
question is a $G$--morphism.---The fixed subspace of $H^{0,1}(L\Bp_1)$ 
can be described more explicitly, for any loop space:

\proclaim{Theorem 0.3}The space $H^{0,1}(L\Bp_1)^G$ is isomorphic to
$C^{k-1}(S^1)^*$, resp.~$W^{k-1,p}(S^1)^*$, if the dual spaces are 
endowed with
the compact--open topology.
\endproclaim

The isomorphisms in Theorem 0.3 are not Diff $S^1$ equivariant.
To remedy this, one is led to introduce the spaces $C^l_r(S^1)$, resp.
$W^{l,p}_r(S^1)$ of differentials $y(t)(dt)^r$ of order $r$ on $S^1$,
of the corresponding regularity; $L_r^p=W_r^{0,p}$.
Then $H^{0,1}(L\Bp_1)^G$ will be Diff $S^1$ equivariantly isomorphic to 
$C_1^{k-1}(S^1)^*$, resp.~$W_1^{k-1,p}(S^1)^*$.

For low regularity loop spaces one can very concretely represent all of 
$H^{0,1}(L\Bp_1)$:

\proclaim{Theorem 0.4} (a) \ If $1\leq p<2$, all of 
$H^{0,1}(L_{1,p}\Bp_1)$ is fixed by
$G$, hence it is isomorphic to $L^{p'}(S^1)$, with $p'=p/(p-1)$.

(b) \ If $1\leq p<\infty$ then $H^{0,1}(L_{1,p}\Bp_1)$ is isomorphic to
$$
\bigoplus_{0\leq n\leq p-1}\frak K_n\otimes 
L_{n+1}^{p/(n+1)}(S^1)^*\approx\bigoplus_{0\leq n\leq p-1}\frak K_n\otimes L_{-n}^{p_n}(S^1),
\qquad p_n={p\over p-1-n},
$$
and so it is the sum of its first $[p]$ isotypical subspaces.
Indeed, the isomorphisms above are $G\times\text{Diff} \ S^1$ 
equivariant, $G$, resp.~Diff $S^1$ acting on one of the factors $\frak K_n,\ 
L_r^q$ naturally, and trivially on the other.
\endproclaim

Again, the dual spaces are endowed with the compact--open topology.

It follows that the infinite dimensional space $H^{0,1}(L_{1,p}\Bp_1)$ 
can
be understood in finite terms, if it is considered as a representation
space of $S^1$. Here $S^1$ acts on itself (by translations), hence also
on $L\Bp_1$ and on $H^{0,1}(L\Bp_1)$. One can read off from
Theorem 0.4 that each irreducible representation of $S^1$ occurs
in $H^{0,1}(L_{1,p}\Bp_1)$ with the same multiplicity $[p]^2$. On the
other hand, for spaces of loops of regularity at least $C^1$, in
$H^{0,1}(L\Bp_1)$ each
irreducible representation of $S^1$ occurs with infinite multiplicity
and, somewhat contrary to earlier expectations, it is not possible
to associate with this cohomology space even a formal character of
$S^1$. This indicates that Dolbeault groups of general loop spaces
$LM$ should be studied as representations of $\text{Diff\,}S^1$ rather
than $S^1$.

The structure of this paper is as follows.
In Sections 1 and 2 we study the space $\frak F$ as a $G$--module.
We connect it with a similar but simpler space of functions that are
required to satisfy only the first two of the three conditions defining
$\frak F$ (Theorem 1.1).
Theorem 1.1 will be needed in proving the isomorphism
$H^{0,1}(L\Bp_1)\approx\Bc\oplus\frak F$, and also in concretely 
representing elements of $\frak F$.
Further, we shall rely on Theorem 1.1 in identifying isotypical 
subspaces of $\frak F$ (Theorems 2.1, 2.2).
This will then prove Theorems 0.2, 0.3, and 0.4, modulo Theorem 0.1.

In Section 3 we introduce a $G$--module $\frak H$ of holomorphic \v 
Cech cocycles of $L\Bp_1$, and prove $H^{0,1}(L\Bp_1)\approx\frak H$ 
(Theorem 3.3).
In Section 4 we construct a morphism $\alpha\colon\frak H\to\frak F$ 
that, in Section 5, is shown to induce an isomorphism $\frak H/\Ker \ 
\alpha\to\frak F$.
Also, dim Ker $\alpha=1$ (Theorem 5.1).
Finally, in Section 6 we show how all this, put together, proves the 
results formulated in this introduction.

\heading 1. The Space $\frak F$\endheading

In this Section and the next we shall study the structure of the space
$\frak F$, independently of any cohomological content.
It will be convenient to allow (but only in this Section!) $k$ to be 
any integer; when $k<0$, elements of $C^k(S^1),\ W^{k,p}(S^1)$ are 
distributions, locally equal to the $-k$'th derivative of functions in 
$C(S^1),\ L^p(S^1)$.
Let $L^-\Bc$ denote the space $C^{k-1}(S^1)$, resp.~$W^{k-1,p} (S^1)$.
We shall write $L^{(-)}\Bc$ to mean either $L\Bc$ or $L^-\Bc$.
Consider the space $\tilde\frak F$ of those $F\in\Cal O(\Bc\times 
L^-\Bc)$ that have properties (1) and (2) of the Introduction.
We shall refer to (2) as additivity.
A function $F\in\Co(\Bc\times L^{(-)}\Bc)$ will be said to be 
posthomogeneous of degree $m$ if $F(\zeta,\cdot)$ is homogeneous of degree $m$ 
for all $\zeta\in\Bc$.
Posthomogeneous degree endows the spaces $\frak F$ and $\tilde\frak F$
with a grading.---All maps below, unless otherwise mentioned, will be
continuous and linear.

\proclaim{Theorem 1.1}The graded linear map
$\tilde\frak F\ni\tilde F\mapsto F\in\frak F$ given by 
$F(\zeta,y)=\tilde F(\zeta,\dot y)$ has a graded right inverse, and its kernel consists 
of functions $\tilde F(\zeta,x)=\const\int_{S^1} x$.
\endproclaim

First we shall consider functions $E\in\frak F$, resp.~$\tilde\frak F$, 
that are
independent of $\zeta$.
We denote the space of these functions $\frak E\subset\Co(L\Bc)$,
resp.~$\tilde\frak E\subset\Co(L^-\Bc)$, graded by degree of 
homogeneity.
Additivity of $E\in\Co(L^{(-)}\Bc)$ implies $E(0)=0$, which in turn 
implies
property (1) of the Introduction. Let
$$
E=\sum^\infty_1 E_m,\qquad E_m(y)=\int_0^1 E(e^{2\pi i\tau} y) 
e^{-2m\pi i\tau}d\tau\tag1.1
$$
be the homogeneous expansion of a general $E\in\Co(L^{(-)}\Bc)$ 
vanishing at 0.
Consider tensor powers $(L^{(-)}\Bc)^{\otimes m}$ of the vector spaces 
$L^{(-)}\Bc$ over $\Bc$.
In particular, $C^\infty(S^1)^{\otimes m}$ is an algebra, and a general 
$(L^{(-)}\Bc)^{\otimes m}$ is a module over it.
Each $E_m$ in (1.1) induces a symmetric linear map
$$
\Ce_m\colon (L^{(-)}\Bc)^{\otimes m}\to\Bc,
$$
called the polarization of $E_m$.
On monomials $\Ce_m$ is defined by
$$
\Ce_m(y_1\otimes\ldots\otimes y_m)={1\over 2^m m!}\sum_{\epsilon_j=\pm 
1}\epsilon_1\ldots\epsilon_m E_m(\epsilon_1 y_1+\ldots+\epsilon_m 
y_m),\tag1.2
$$
and then extended by linearity.
Thus $E_m(y)=\Ce_m(y^{\otimes m})$.---We shall call $w\in 
(L^{(-)}\Bc)^{\otimes m}$ degenerate if it is a linear combination of monomials 
$y_1\otimes\ldots\otimes y_m$ with one $y_j=1$.

\proclaim{Lemma 1.2}(a)\, $E$ is additive if and only if 
$\Ce_m(y_1\otimes\ldots\otimes y_m)=0$ whenever
$\bigcap_1^m\supp \, y_j~=~\emptyset$.

(b) $E(y+\text{const})=E(y)$ if and only if $\Ce_m(w)=0$ whenever $w$ 
is degenerate.

\endproclaim

\demo{Proof}(a) Clearly $E$ is additive precisely when all the $E_m$ 
are, whence it suffices to prove the claim when $E$ itself is 
homogeneous, of degree $m$, say.
In this case $\Ce_n=0,\ n\not= m$.
Denoting $\Ce_m$ by $\Ce$, it is also clear that the condition on $\Ce$ 
implies $E$ is additive.
We show the converse by induction on $m$, the case $m=1$ being obvious.
Let $x,y\in L^{(-)}\Bc$ have disjoint support, so that
$$
\Ce((x+y)^{\otimes m})=\Ce(x^{\otimes m})+\Ce(y^{\otimes m}).\tag1.3
$$
Write $\lambda x$ for $x$ and separate terms of different degrees in 
$\lambda$ to find $\Ce(x\otimes\ldots\otimes y)=0$, which settles the 
case $m=2$.
Now if the case $m-1\geq 2$ is already settled, take a $z\in 
L^{(-)}\Bc$ with supp $y\cap\supp\,z=\emptyset$, and write $x+\lambda z$ for $x$ 
in (1.3).
Considering the terms linear in $\lambda$ we obtain
$$
\Ce(z\otimes(x+y)^{\otimes m-1})=\Ce(z\otimes x^{m-1})+\Ce(z\otimes 
y^{m-1}),\tag1.4
$$
the last term being 0.
The same will hold if supp $x\cap\supp \, z=\emptyset$.
Since any $z\in L^{(-)}\Bc$ can be written $z'+z''$ with the support of 
$z'$ (resp.~$z''$) disjoint from the support of $x$ (resp.~$y$), (1.4) 
in fact holds for all $z$.
By the induction hypothesis applied to $\Ce(z\otimes\cdot)$
$$
\Ce(z\otimes y_2\otimes\ldots\otimes y_m)=0,\qquad\text{if}\quad
\bigcap^m_2\supp \ y_j=\emptyset.
$$
Suppose now $\bigcap^m_1\supp \, y_j=\emptyset$ and write $y_1=y'+y''$ 
with $y'=0$ near $\bigcap_{j\not= 2} \supp\,y_j$ and $y''=0$ near 
$\bigcap_{j\not= 3}\supp\, y_j$.
Then
$$
\Ce(y_1\otimes\ldots\otimes y_m)=\Ce(y'\otimes\ldots\otimes 
y_m)+\Ce(y''\otimes\ldots\otimes y_m)=0.
$$

(b) Again we assume $E$ is $m$--homogeneous, and again one implication 
is trivial.
So assume $\Ce((y+1)^{\otimes m})=\Ce(y^{\otimes m})$, where 
$\Ce=\Ce_m$.
Differentiating both sides in the directions $y_2,\ldots,y_m$ and 
setting $y=0$ we obtain $\Ce(1\otimes y_2\otimes\ldots\otimes y_m)=0$, 
whence the claim follows.
\enddemo

\proclaim{Proposition 1.3}The graded map $\tilde\frak E\ni\tilde 
E\mapsto E\in\frak E$ given by $E(y)=\tilde E(\dot y)$ has a graded right 
inverse, and its kernel is spanned by $\tilde E(x)=\int_{S^1} x$.
\endproclaim

We shall write $\int x$ for $\int_{S^1} x$.
\demo{Proof}(a) To identify the kernel, because of homogeneous 
expansions, it will suffice to deal with homogeneous $\tilde E$.
So assume $\tilde E\in\tilde\frak E$ is homogeneous of degree $m$ and 
$\tilde E(\dot y)=0$ for all $y\in L\Bc$.
Its polarization $\tilde\Ce$ satisfies
$\tilde\Ce(\dot y_1\otimes\ldots\otimes\dot y_m)=0$.
If $m=1$, this implies $\tilde E(x)=$ const $\int x$, so from now on we 
assume $m\geq 2$, and first we prove by induction that 
$\tilde\Ce(x_1\otimes\ldots\otimes x_m)=\const\,\prod\int x_j$.
Suppose we already know this for $m-1$.
Then
$$
\tilde\Ce(\dot y \otimes x_2\otimes\ldots\otimes x_m)=c(\dot 
y)\prod^m_2\int x_j.
$$
With arbitrary $x_1\in L^-\Bc$ the function $x_1-\int x_1$ is of form 
$\dot y$, so $x_1=\dot y+\int x_1$, and
$$
\tilde\Ce(x_1\otimes\ldots\otimes x_m)=l(x_1)\prod_2^m\int 
x_j+\tilde\Ce(1\otimes x_2\otimes\ldots\otimes x_m) \int x_1,\tag1.5
$$
where $l(x_1)=c(x_1-\int x_1)$ is linear in $x_1$.
If $\int x_1=0$ and $\supp\,x_1\not= S^1$, then we can choose 
$x_2,\ldots$ so that $\bigcap_1^m\supp \, x_j=\emptyset$ but $\int x_j\not= 0,\ 
j\geq 2$.
This makes the left hand side of (1.5) vanish by Lemma 1.2a, and gives 
$l(x_1)=0$.
Since any $x_1\in L^-\Bc$ with $\int x_1=0$ can be written $x_1=x'+x''$ 
with $\int x'=\int x''=0$ and $\supp \,x',\supp \,x''\not= S^1$, it 
follows that $l(x_1)=0$ whenever $\int x_1=0$.
Hence $l(x_1)=\const\int x_1$.
In particular, the first term on the right of (1.5) is symmetric in 
$x_j$.
Therefore the second term must be symmetric too, which implies this 
term is $\const\prod_1^m\int x_j$.
Thus $\tilde E(x)=\const(\int x)^m$.

Yet for $m\geq 2$ $\tilde E(x)=\const(\int x)^m$ is additive only if it 
is identically zero; so that indeed $\tilde E(x)=\const\int x$, as 
claimed.

(b) To construct the right inverse, consider $E\in\frak E$ with 
homogeneous expansion (1.1).
We shall construct $m$--homogeneous polynomials $\tilde 
E_m\in\tilde\frak E$ such that $E_m(y)=\tilde E_m(\dot y)$; the case $m=1$ being 
obvious, we assume $m\geq 2$.
Let us say that an $n$--tuple of functions $\rho_\nu\colon S^1\to\Bc$ 
is centered if $\bigcap_1^n\supp\,\rho_\nu\not=\emptyset$.
We start by fixing a $C^\infty$ partition of unity $\sum_{\rho\in 
P}\rho=1$ on $S^1$ such that each $\supp\,\rho$ is an arc of length $<1/4$. 
This
implies that $\bigcup_1^n\supp \,\rho_\nu$ is an arc of length $<1/2$ 
if
$\rho_1,\ldots,\rho_n\in P$ are centered.
Given $x\in L^-\Bc$, for each centered $R=(\rho_1,\ldots,\rho_n)$ in 
$P$ construct $y_R\in L\Bc$ so that $\dot y_R=x$ on a neighborhood of
$\bigcup_1^n\supp\,\rho_\nu$, making sure that $y_R=y_Q$ if $Q$ and $R$
agree as sets.
For noncentered $n$--tuples $R$ in $P$ let $y_R\in L\Bc$ be arbitrary.
We shall refer to the $y_R$ as local integrals.

If $Q,R$ are centered tuples in $P$ then
$$
y_Q-y_R=c_{QR}=\text{constant}\qquad\text{on}\quad
(\bigcup_{\rho\in Q}\supp \ \rho)\cap (\bigcup_{\rho\in R}\supp \ 
\rho).\tag1.6
$$
When the intersection in (1.6) is empty, or $Q$ or $R$ are not 
centered, fix $c_{QR}\in\Bc$ arbitrarily.
Define
$$
v_{QR}=m\int^{c_{Q R}}_0(y_R+\tau)^{\otimes m-1} d\tau\in 
(L\Bc)^{\otimes m-1},\tag1.7
$$
and with the polarization $\Ce_m$ of $E_m$ from (1.2) consider
$$
\Ce_m\bigg(\sum_{R=(\rho_1,\ldots,\rho_m)}(\rho_1\otimes\ldots\otimes\rho_m)
\bigl(y_R^{\otimes m}+1\otimes\sum_{S=(\sigma_2,\ldots,\sigma_m)}
(\sigma_2\otimes\ldots\otimes\sigma_m)v_{SR}\bigr)\bigg);\tag1.8
$$
we sum over all $m$--tuples $R$ and $(m-1)$--tuples $S$ in $P$.
(We will not need it, but here is an explanation of (1.8).
Say that tensors $w,w'\in L^{(-)}\Bc^{\otimes m}$ are congruent, 
$w\equiv w'$, if $w-w'$ is the sum of a degenerate tensor and of monomials 
$x_1\otimes\ldots\otimes x_m$ with $\bigcap\supp \ x_j=\emptyset$.
Denote by $\partial^m$ the linear map $(L\Bc)^{\otimes m}\to 
(L^-\Bc)^{\otimes m}$ defined by $\partial^m(y_1\otimes\ldots\otimes y_m)=\dot 
y_1\otimes\ldots\otimes\dot y_m$.
Then the symmetrization of the argument of $\Ce_m$ in (1.8) is a 
solution $w$ of the congruence $\partial^m w\equiv x^{\otimes m}$, in fact it 
is the unique symmetric solution, up to congruence.
It follows that for the $\tilde E_m$ sought, $\tilde E_m(x)$ must be 
equal to $\Ce_m(w)$, which, in turn, equals (1.8).)

We claim that the value in (1.8) depends only on $x$ (and $\Ce_m$), but 
not on the partition of unity $P$ and the local integrals  $y_R$.
Indeed, suppose first that the local integrals $y_R$ are changed to 
$\hat y_R$, so that the $c_{QR}$ change to $\hat c_{QR}$ and $v_{QR}$ to 
$\hat v_{QR}$; but we do not change $P$.
There are $c_R\in\Bc$ such that for all centered $R$
$$
\hat y_R=y_R+c_R\qquad\text{ on }\bigcup_{\rho\in R}\supp \ \rho.
$$
Let
$$
u_R=m\int_0^{c_R}(y_R+\tau)^{\otimes m-1}d\tau.\tag1.9
$$
Clearly $\hat c_{QR}=c_{QR}+c_Q-c_R$ if $Q\cup R$ is centered.
In this case one computes also
$$
\alignedat 2
\frac{1}{m}\hat v_{QR}&=\int_0^{\hat c_{QR}}
&&(\hat y_R+\tau)^{\otimes m-1}
d\tau \\
&=\int_0^{c_{QR}}&&(\hat y_R-c_R+\tau)^{\otimes m-1}d\tau-
\int_0^{c_R}(\hat y_R-c_R+\tau)^{\otimes m-1}d\tau\\
& \ &&+
\int_0^{c_Q}(\hat y_R-c_R+c_{QR}+\tau)^{\otimes m-1}d\tau.
\endalignedat\tag1.10
$$

Because of Lemma 1.2a, in (1.8) only centered $R$, and such $S$ that 
$R\cup S$ is centered, will contribute.
When $y_R^{\otimes m}$ is changed to $\hat y_R^{\otimes m}$, the 
corresponding contributions change by
$$
\align
\sum_R\Ce_m\bigg(\int_0^{c_R}(\rho_1 & \otimes\ldots\otimes\rho_m){d\over d\tau}(y_R+\tau)^{\otimes m}dt\bigg) \\
&=\sum_R\Ce_m\bigg(
\int_0^{c_R}(\rho_1\otimes\ldots\otimes\rho_m)(m\otimes(y_R+\tau)^{\otimes 
m-1})d\tau\bigg) \\
&=\sum_R\Ce_m((\rho_1\otimes\ldots\otimes\rho_m)(1\otimes u_R)).
\endalign
$$
When $v_{QR}$ is changed to $\hat v_{QR}$, in view of (1.10), (1.6), 
and (1.9), the contribution of the terms in the double sum in (1.8) 
changes by
$$
\gathered
\Ce_m\bigg((m\rho_1\otimes\rho_2\sigma_2\otimes\ldots
\otimes\rho_m\sigma_m)
\bigl(\int_0^{c_S}(y_S+\tau)^{\otimes m-1}d\tau-
\int_0^{c_R}(y_R+\tau)^{\otimes m-1 }d\tau\bigr)\bigg)\\
=\Ce_m((\rho_1\otimes\rho_2\sigma_2\otimes\ldots\otimes\rho_m\sigma_m)(1\otimes 
u_{S }-1\otimes u_R)).
\endgathered
$$
The net change in (1.8) is therefore
$$
\aligned
\Ce_m\bigg(\sum_{R,S}(\rho_1&
\otimes\rho_2\sigma_2\otimes\ldots\otimes\rho_m\sigma_m)(1\otimes 
u_{S})\bigg)=\\
&\Ce_m\bigg(\sum_S(1\otimes\sigma_2\otimes\ldots\otimes\sigma_m)
(1\otimes u_{S })\bigg)=0
\endaligned
$$
by Lemma 1.2b, as needed.

Now to pass from $P$ to another partition of unity $P'$, introduce 
$\Pi=\{\rho\rho'\colon p\in P,\rho'\in P'\}$.
One easily shows that $P$ and $\Pi$ give rise to the same value in 
(1.8), hence so do $P$ and $P'$.
Therefore (1.8) indeed depends only on $x$, and we define $\tilde 
E_m(x)$ by this value.
We proceed to check that $\tilde E_m$ has the required properties.

If $x=\dot y$ then all $y_R$ can be chosen $y$, and (1.8) gives $\tilde 
E_m(\dot y)=E_m(y)$.
Next suppose $x',x''\in L^-\Bc$ have disjoint support, and $x=x'+x''$.
If the supports of all $\rho\in P$ are sufficiently small, then the 
local integrals $y'_R,\ y''_R$ of $x', x''$ can be chosen so that for each 
$R$ one of them is 0.
Hence the local integrals $y_R=y'_R+y''_R$ of $x$ will satisfy
$y_R^{\otimes m}=y_R^{\prime\otimes m} +y_R^{\prime\prime\otimes m}$, 
whence $\tilde E_m(x)=\tilde E_m(x')+\tilde E_m(x'')$ follows.

To show that $\sum\tilde E_m$ is convergent and represents a 
holomorphic function, note that $\tilde E_m(x)$ is the sum of terms
$$
\gathered
\Ce_m(\rho_1 y_R\otimes\ldots\otimes\rho_m y_R)\qquad\text{and}\\
\int_0^1\Ce_m(\rho_1 c_{SR}\otimes\rho_2\sigma_2(y_R+
c_{SR}\tau)\otimes\ldots\otimes\rho_m\sigma_m(y_R+
c_{SR}\tau))d\tau
\endgathered\tag1.11
$$
(we have substituted $c_{QR}\tau$ for $\tau$ in (1.7)).
Since $y_R\in L\Bc$ and $c_{QR}\in\Bc$ can be chosen to depend on $x$ 
in a continuous linear way, each $\tilde E_m$ is a homogeneous 
polynomial of degree $m$.
Furthermore, let $K\subset L^-\Bc$ be compact.
For each $x\in K,\ m\in\Bn$, and $m$--tuples $Q,R$ in $P$ we can choose 
$y_R$ and $c_{QR}$ so that all the functions
$$
\rho c_{QR},\ \rho\rho'(y_R+c_{QR}\tau),
$$
$\rho,\rho'\in P,\ 0\leq \tau\leq 1$, belong to some compact $H\subset 
L\Bc$.
By passing to the balanced hull, it can be assumed $H$ is balanced.
If $\lambda > 0$, (1.1) implies
$$
\max_H |E_m|=\lambda^{-m}\max_{\lambda H} 
|E_m|\leq\lambda^{-m}\max_{\lambda H}|E|=A\lambda^{-m},
$$
so that by (1.2)
$$
|\Ce_m(z_1\otimes\ldots\otimes z_m)|\leq A(m^m/m!)\ \lambda^{-m}\leq 
A(e/\lambda)^m,
$$
if each $z_\mu\in H$.
Thus each term in (1.11) satisfies this estimate.
If $|P|$ denotes the cardinality of $P$, we obtain, in view of (1.8)
$$
\max_K|\tilde E_m|\leq (|P|^m+m|P|^{2m-1})A(e/\lambda)^m.
$$
Choosing $|\lambda|> e|P|^2$ we conclude that $\sum\tilde E_m$ 
uniformly converges on $K$, and, $K$ being arbitrary, $\tilde E=\sum\tilde E_m$ 
is holomorphic.
By what we have already proved for $\tilde E_m,\ \tilde E\in\tilde\frak 
E$, and $\tilde E(\dot y)=E(y)$.
The above estimates also show that the map $E\to\tilde E$ is continuous 
and linear, which completes the proof of Proposition 1.3.
\enddemo

Now consider an $F\in\Co(\Bc\times L^{(-)}\Bc)$ and its posthomogeneous 
expansion
$$
F=\sum^\infty_0 F_m,\quad F_m(\zeta,y)=\int_0^1 F(\zeta,e^{2\pi i\tau} 
y)e^{-2m\pi i\tau}d\tau.\tag1.12
$$

\proclaim{Proposition 1.4}The function $F$ satisfies condition (1) of 
the Introduction if and only if each $F_m$ is a polynomial in $\zeta$, 
of degree $\leq 2m-2$ (in particular, $F_0=0$).
\endproclaim

\demo{Proof}As $F$ satisfies (1) precisely when each $F_m$ does, the 
statement is obvious.
\enddemo

\demo{Proof of Theorem 1.1}Apply Proposition 1.3 on each slice 
$\{\zeta\}\times L^{(-)}\Bc$. Accordingly, an $\tilde F$ in the kernel is 
posthomogeneous
of degree 1, hence, by Proposition 1.4, independent of $\zeta$. Thus
indeed $\tilde F(\zeta,x)=\const\int x$.
Further, the slicewise right inverse applied to $F\in\frak F$ clearly 
produces an additive $\tilde F\in\Co(\Bc\times L\Bc)$.
To see that $\tilde F$ also verifies condition (1) of the Introduction, 
expand $F$ in a posthomogeneous series
$$
F(\zeta,y)=\sum^\infty_{m=1}F_m(\zeta,y)=\sum^\infty_{m=1}\sum^{2m-2}_{\nu=0}\zeta^\nu 
E_{m\nu}(y),\tag1.13
$$
by Proposition 1.4, so that
$$
\tilde F(\zeta,x)=\sum^\infty_{m=1}\sum^{2m-2}_{\nu=0}\zeta^\nu\tilde 
E_{m\nu}(x),
$$
with $\tilde E_{m\nu}$ $m$--homogeneous.
Again by Proposition 1.4, $\tilde F$ verifies condition (1), and so
$\in\tilde\frak F$.
\enddemo

Theorem 1.1 can be used effectively to describe elements of the space 
$\frak F$.
With ulterior motives we switch notation $m=n+1$, and consider a 
homogeneous polynomial $\tilde E\in\Co(L^-\Bc)$ of degree $n+1\geq 1$.
Its polarization $\Ce$ defines a distribution $D$ on the torus 
$(S^1)^{n+1}=T$.
Indeed, denote the coordinates on $T$ by $t_j\in\Bbb R/\Bbb Z$ and set
$$
\langle D,\ \prod_{j=0}^n e^{2\pi i\nu_j 
t_j}\rangle=\Ce(x_0\otimes\ldots\otimes x_n),\qquad x_j(\tau)=e^{2\pi i\nu_j\tau},\ 
\nu_j\in\Bz.\tag1.14
$$
Since $\tilde E$ is continuous,
$$
|\Ce(x_0\otimes\ldots\otimes x_n)|\leq 
c\prod\|x_j\|_{C^q(S^1)}\qquad\text{ with
some }c>0  \text{ and }q\in\Bn.
$$
Hence (1.14) can be estimated, in absolute value, by 
$c'\prod_j(1+|\nu_j|)^q$, and it follows by Fourier expansion that $D$ extends to a 
unique linear form on $C^\infty(T)$.
Clearly, $D$ is symmetric, i.e., invariant under permutation of the 
factors $S^1$ of $T$.
Also,
$$
\Ce(x_0\otimes\ldots\otimes x_n)=\langle D,x_0\otimes\ldots\otimes 
x_n\rangle,\tag1.15
$$
if on the right $x_0\otimes\ldots\otimes x_n$ is identified with the 
function $\prod x_j(t_j)$.

Assume now $\tilde E\in\tilde\frak E$.
Lemma 1.2a implies $D$ is supported on the diagonal of $T$.
The form of distributions supported on submanifolds is in general well 
understood; in the case at hand, e.g.~[H, Theorem 2.3.5] gives that $D$ 
is a
finite sum of distributions of form
$$
C^\infty(T)\ni\rho\mapsto\biggl\langle\Psi,\ 
{{\partial^{\alpha_1+\ldots+\alpha_n}\rho}\over{\partial t_1^{\alpha_1}\ldots\partial 
t_n^{\alpha_n}}} \bigg|_{\text{diag}}\biggr\rangle,\quad \alpha_j\geq 0,
$$
where $\Psi$ is a distribution on the diagonal of $T$.
In view of Theorem 1.1 and (1.12)--(1.13) we therefore proved

\proclaim{Theorem 1.5}The restriction of an $(n+1)$--posthomogeneous 
$F\in\frak F$, resp.~$\tilde\frak F$, to $\Bc\times C^\infty(S^1)$ is a 
finite sum of functions of form
$$
f(\zeta,y)=\zeta^\nu\biggl\langle\Phi,\prod^n_0 
y^{(d_j)}\bigg\rangle,\qquad \nu\leq 2n,\ d_j\geq d_0=1,\text{ resp. }0,
$$
where $\Phi$ is a distribution on $S^1$.
For a general $F\in\frak F$, resp.~$\tilde\frak F$, the restriction 
$F|\Bc\times C^\infty(S^1)$ is the limit, in the topology of
$\Cal O(\Bbb C\times C^\infty(S^1))$, of finite sums of the above 
functions.
\endproclaim

\heading 2.\ The $G$--action on $\frak F$\endheading

For $g\in G$ let $J_g(\zeta)=dg\zeta/d\zeta$.
By considering the posthomogeneous expansion (1.12)--(1.13) of 
$F\in\frak F$,
resp.~$\tilde\frak F$, one checks that the function $gF$ defined by
$$
(gF)(\zeta,y)=F(g\zeta,\ y/J_g(\zeta))J_g(\zeta)\tag2.1
$$
extends to all of $\Bc\times L^{(-)}\Bc$, and the extension (also 
denoted
$gF$) belongs to $\frak F$, resp.~$\tilde\frak F$.
The action thus defined makes $\frak F,\tilde\frak F$ holomorphic 
$G$--modules.
The $n$'th isotypical subspace $\frak F^n$, resp.~$\tilde\frak F^n$ is 
the subspace of $(n+1)$--posthomogeneous functions.
In this section we shall describe the space $\frak F^0$, and, for 
$W^{1,p}$ loop spaces, the spaces $\frak F^n$ as well, $n\geq 1$.

\proclaim{Theorem 2.1}$\frak F^0\approx (L^-\Bc)^*/\Bc$, the dual 
endowed with the compact--open topology.
If $L^-\Bc$ is interpreted as the space of one--forms on $S^1$ of the 
corresponding regularity, then the isomorphism is Diff $S^1$ 
equivariant.
\endproclaim

\demo{Proof}Indeed, the map $(L^-\Bc)^*=\tilde\frak F^0\to\frak F^0$ 
associating with $\Phi\in (L^-\Bc)^*$ the function $F(y)=\langle\Phi,\dot 
y\rangle$ (or $\langle\Phi,dy\rangle$) has one dimensional kernel and a 
right inverse by Theorem 1.1.
\enddemo

\proclaim{Theorem 2.2}In the case of $W^{1,p}$ loop spaces $\frak 
F=\bigoplus_{n\leq p-1}\frak F^n$.
Furthermore
$$
\frak K_n\otimes L^{p/(n+1)}(S^1)^*\approx \frak F^n,\qquad 1\leq n\leq 
p-1,
$$
as $G$--modules, $G$ acting on $L^{p/(n+1)}(S^1)^*$ trivially.
Indeed, the map $\varphi\otimes\Phi\mapsto F$ given by
$$
F(\zeta,y)=\psi(\zeta)\langle\Phi,\dot 
y^{n+1}\rangle,\qquad\varphi(\zeta)=\psi(\zeta)(d\zeta)^{-n},\tag2.2
$$
induces the isomorphism above.
(To achieve Diff $S^1$ equivariant isomorphism, replace 
$L^{p/(n+1)}(S^1)$ by the space $L_{n+1}^{p/(n+1)}(S^1)$ of $(n+1)$--differentials.)
\endproclaim
We shall need a few auxiliary results to prove the theorem.
\proclaim{Lemma 2.3}Let $m\geq 2$ be an integer and $\Psi$ a 
distribution on $S^1$.
If the function
$$
C^\infty(S^1)\ni x\mapsto \langle\Psi,x^m\rangle\in\Bc\tag2.3
$$
extends to a homogeneous polynomial $E$ on $L^p(S^1)$ then $\Psi\equiv 
0$, or $m\leq p$ and $\Psi$ extends to a form $\Phi$ on $L^{p/m}(S^1)$.
In the latter case the map $E\mapsto\Phi$ is continuous 
linear.\endproclaim

\demo{Proof}There is a constant $C$ such that
$$
|\langle \Psi,x^m\rangle|=|E(x)|\leq C(\int |x|^p)^{m/p},\qquad x\in 
C^\infty(S^1).\tag2.4
$$
Let $z\in C^\infty(S^1)$ be real valued and 
$x_\epsilon=(z+i\epsilon)^{1/m}\in C^\infty(S^1)$.
By (2.4)
$$
|\langle\Psi,z\rangle|=\lim_{\epsilon\to 0} 
|\langle\Psi,x_\epsilon^m\rangle|\leq C(\int |z|^{p/m})^{m/p}.
$$
As the same estimate holds for imaginary $z$, it will hold for a 
general $z\in C^\infty(S^1)$ too, perhaps with a different $C$.
Therefore $\Psi$ extends to a form $\Phi$ on $L^{p/m}(S^1)$.
Unless $p\geq m$, $\Phi=0$ by Day's theorem [D].
With $z\in L^{p/m}(S^1)$, any choice of measurable $m$'th root 
$z^{1/m}$, and
$y_\varepsilon\in C^\infty(S^1)$ converging to $z^{1/m}$ in $L^p$,
$$
\langle\Phi,z\rangle=\lim_{\varepsilon\to 
0}\langle\Phi,y_\varepsilon^m\rangle=\lim_{\varepsilon\to 0} E(y_\varepsilon)=E(z^{1/m}).
$$
This shows that $\Phi$ is uniquely determined by $E$, and depends
continuously and linearly on $E$.
\enddemo

In the rest of this section we work with $W^{1,p}$ loop spaces.
Write $\frak E^n\subset \frak E,\ \tilde\frak E^n\subset\tilde\frak E$ 
for the space of $(n+1)$--homogeneous functions.

\proclaim{Lemma 2.4}If $m\geq 2$ and $E\in\tilde\frak 
E^{m-1}\subset\Cal O(L^p(S^1))$, then $E(x)=\langle\Phi,x^m\rangle$ with a unique 
$\Phi\in L^{p/m}(S^1)^*$.
In particular, $E=0$ if $m>p$.
Also, the map $E\mapsto \Phi$ is an isomorphism between $\tilde\frak 
E^{m-1}$ and $L^{p/m}(S^1)^*$.
\endproclaim

\demo{Proof}We shall prove by induction, first assuming $m=2$.
By Theorem 1.5 there are distributions $\Phi_\alpha$ so that
$$
E(x)=\sum^d_{\alpha=0}\langle\Phi_\alpha,xx^{(\alpha)}\rangle,\quad 
x\in C^\infty(S^1).
$$
Now any $x^{(\alpha)}x^{(\beta)}$ will be a linear combination of 
expressions $(x^{(j)} x^{(j)})^{(h)}$, as one easily proves by induction of 
$|\alpha-\beta|$.
It follows that $E$ can be written with distributions $\Psi_j$ as
$$
E(x)=\sum^d_{j=0}\langle\Psi_j,(x^{(j)})^2\rangle,\quad x\in 
C^\infty(S^1).\tag2.5
$$
Next we show that $d=0$.

Indeed, assuming $d>0$, for fixed $x\in C^\infty(S^1)$
$$
E(\cos\lambda x)+E(\sin\lambda x)=\lambda^{2d}\langle\Psi_d,\dot 
x^{2d}\rangle+\sum^{2d-1}_{j=0}c_j(x)\lambda^d\tag2.6
$$
is a polynomial in $\lambda$.
For fixed $\lambda\in\Bc$ the maps $x\mapsto\cos\lambda x,\ 
x\mapsto\sin\lambda x$ map the Banach algebra $W^{1,1}(S^1)$ holomorphically into 
itself, hence into $L^p(S^1)$.
Therefore the left hand side of (2.6) extends to $W^{1,1}(S^1)$, and 
$\langle\Psi_d,\dot x^{2d}\rangle$ must also.
The extension of this latter will be an additive, $2d$--homogeneous 
polynomial $E'$ on $W^{1,1}(S^1)$, satisfying $E'(x+\const)=E'(x)$.
By Proposition 1.3 there is therefore a unique additive 
$2d$--homogeneous polynomial $\tilde E$ on $W^{0,1}(S^1)=L^1(S^1)$ such that 
$E'(x)=\tilde E(\dot x)$.
Since the restriction $\tilde E|C^\infty(S^1)$ is also unique,
$$
\tilde E(x)=\langle\Psi_d, x^{2d}\rangle,\qquad x\in C^\infty(S^1).
$$
In particular, the expression on the right continuously extends to 
$L^1(S^1)$.
By virtue of Lemma 2.3, $\Psi_d\equiv 0$.
Thus (2.5) reduces to $E(x)=\langle\Psi,x^2\rangle,\ x\in 
C^\infty(S^1)$, and by another application of Lemma 2.3, $\Psi$ extends to a form 
$\Phi$ on $L^{p/2}(S^1)$.

Now assume the Lemma is known for degree $m-1\geq 2$, and consider an 
$E\in\tilde\frak E^{m-1}$ and its polarization $\Cal E$.
For fixed $x_1\in C^\infty(S^1)$ the inductive assumption implies that 
there is a distribution $\Theta$ such that $\Cal 
E(x_1\otimes\ldots\otimes x_m)=\langle\Theta,\prod^m_2 x_j\rangle$; in particular,
$$
\Cal E(x_1\otimes\ldots\otimes x_m)=\Cal E(x_1\otimes\prod^m_2 
x_j\otimes 1\otimes\ldots\otimes 1),\qquad x\in C^\infty(S^1).
$$
The case $m=2$ now gives a distribution $\Psi$ such that $\Cal 
E(x_1\otimes\ldots\otimes x_m)=\langle\Psi,\prod^m_1 x_j\rangle$.
We conclude by Lemma 2.3:\ $\Psi$ extends to $\Phi\in L^{p/m}(S^1)^*$, 
and $\Phi=0$ unless $m\leq p$.
It is clear that $\Phi$ is uniquely determined by $E$, and the map 
$\tilde\frak E^{m-1}\ni E\mapsto \Phi\in L^{p/m}(S^1)^*$ is an isomorphism.
\enddemo

\demo{Proof of Theorem 2.2}To construct the inverse of the map defined 
by (2.2), write an arbitrary $F\in\frak F^n,\ n\geq 1$, as
$$
F(\zeta,y)=\sum^{2n}_{\nu=0}\zeta^\nu E_\nu(y),\quad E_\nu\in\frak E^n,
$$
cf. Proposition 1.4,
and find the unique $\tilde E_\nu\in\tilde\frak E^n$ so that 
$E_\nu(y)=\tilde E_\nu(\dot y)$, see Proposition 1.3.
By Lemma 2.4 there are unique $\Phi_\nu\in L^{p/(n+1)}(S^1)^*$ such 
that $\tilde E_\nu(x)=\langle\Phi_\nu,x^{n+1}\rangle$.
If $p<n+1$ then $\Phi_\nu=0$ and so $\frak F^n=(0)$.
Otherwise the map
$$
\frak F^n \ni 
F\mapsto\sum^{2n}_0\zeta^\nu(d\zeta)^{-n}\otimes\Phi_\nu\in\frak K_n\otimes L^{p/(n+1)}(S^1)^*
$$
is the inverse of the map given in (2.2), so (2.2) indeed induces an 
isomorphism.
Finally, the posthomogeneous expansion of an arbitrary $F\in\frak F$ is
$$
F=\sum^\infty_0 F_n=\sum^{[p-1]}_0 F_n,
$$
which completes the proof.
\enddemo

\heading{3. Cuspidal cocycles}\endheading

In this section we shall construct an isomorphism between 
$H^{0,1}(L\Bp_1)$ and a space of holomorphic \v Cech cocycles on $L\Bp_1$.
We represent $\Bp_1$ as $\Bc\cup\{\infty\}$.
Constant loops constitute a submanifold of $L\Bp_1$, that we identify 
with $\Bp_1$.
If $a,b,\ldots\in\Bp_1$, set $U_{ab\ldots}=\Bp_1\backslash 
\{a,b,\ldots\}$.
Thus $LU_a,\ a\in\Bp_1$, form an open cover of $L\Bp_1$, with 
$LU_\infty=L\Bc$ a Fr\'echet algebra.
If $g\in G$ then $g(LU_a)=LU_{ga}$.

Suppose we are given $v\colon\Bp_1\to\Bc$, finitely many 
$a,b,\ldots\in\Bp_1$,
and a function $u:LU_{ab\ldots}\to\Bc$. If $\infty$ is among 
$a,b,\ldots$,
let us say that $u$ is $v$-cuspidal at $\infty$ if $u(x+\lambda)\to 
v(\infty)$ as $\Bc\ni\lambda\to\infty$, for all $x\in LU_{ab\ldots}$; and 
in
general, that $u$ is $v$-cuspidal if $g^* u$ is $g^*v$-cuspidal at
$\infty$ for all $g\in G$ that maps one of $a,b,\ldots$ to $\infty$.
When $v\equiv 0$ we simply speak of cuspidal functions.

\proclaim{Proposition 3.1}Given a closed $f\in C_{0,1}^\infty(L\Bp_1)$ 
and
$v\in C^\infty(\Bp_1)$ such that $\overline\partial v=f|\Bp_1$, for 
each
$a\in\Bp_1$ there is a unique $v$-cuspidal $u_a\in C^\infty (LU_a)$ 
that
solves $\overline\partial u_a=f|LU_a$.
Furthermore, $u_a|U_a=v|U_a$, and $u(a,x)=u_a(x)$ is smooth in $(a,x)$,
holomorphic in $a$.
\endproclaim

\demo{Proof}Uniqueness follows since for fixed $g\in G,\ y\in L\Bc$, on 
the line $\{g(y+\lambda)\colon\lambda\in\Bp_1\}$ the 
$\overline\partial$ equation is uniquely solvable up to an additive constant, which 
constant is determined by the cuspidal condition.
To construct $u_a$, fix a $g\in G$ with $g\infty=a$, let $Y=\{y\in 
L\Bc\,
\colon y(0)=0\}$ and
$$
P_g\colon\Bp_1\times Y\ni (\lambda,y)\mapsto g(y+\lambda)\in L\Bp_1,
$$
a biholomorphism between $\Bc\times Y$ and $LU_a$.
Setting $f_g=P_g^* f$, by [L1, Theorem 5.4] on the $\Bp_1$ bundle 
$\Bp_1\times Y$ the equation $\overline\partial u_g=f_g$ has a unique smooth 
solution satisfying $u_g(\infty,x)=v(a)$.
It follows that $u_a=(P_g^{-1})^* (u_g|\Bc\times Y)$ solves 
$\overline\partial u_a=f|LU_a$.
Also, $g^* u_a$ is $g^* v$--cuspidal at $\infty$.
On $U_a$ both $u_a$ and $v$ solve the same 
$\overline\partial$--equation,
and have the same  limit at $a$, hence $u_a|U_a=v|U_a$.

One can also consider
$$
P\colon\Bp_1\times G\times Y\ni(\lambda,g,y)\mapsto g(y+\lambda)\in 
L\Bp_1
$$
and $f'=P^*f$.
Again by [L1, Theorem 5.4], on the $\Bp_1$ bundle $\Bp_1\times G\times 
Y$ the equation $\overline\partial u'=f'$ has a smooth solution 
satisfying $u'(\infty,g,x)=v(g\infty)$.
Uniqueness of $u_g$ implies $u'(\lambda,g,x)=u_g(\lambda,x)$, whence 
$u_g(\lambda,x)$ depends smoothly on $(\lambda,g,x)$, and $u_a(x)$ on 
$(a,x)$.
Furthermore, $u'$ is holomorphic on $P^{-1}(x)$ for any $x$.
In particular, if $g\in G$ with $g\infty=a$ is chosen to depend 
holomorphically on $a$ (which can be done locally), then it follows that 
$u_a(x)=u'(g^{-1}x(0), g,g^{-1}x-g^{-1}x(0))$ is holomorphic in $a$.
\enddemo

Since $f$ determines
$v$ up to an additive constant, we can uniquely associate
with $f$ the \v Cech cocycle $\frak f=(u_a-u_b: a,b\in\Bp_1)$.  The
components of $\frak f$ are cuspidal holomorphic functions on 
$LU_{ab}$. One
easily verifies

\proclaim{Proposition 3.2} $f$ is exact if and only if $\frak f=0$. 
Hence
$\frak f$ depends only on the cohomology class $[f]\in 
H^{0,1}(L\Bp_1)$.
The components $h_{ab}([f],x)$ of $\frak f$ depend holomorphically on 
$a,b\in\Bp_1$ and $x\in LU_{ab}$, and satisfy the transformation formula
$$
h_{ga,gb}([f],gx)=h_{ab}(g^*[f],x),\qquad g\in G,\ x\in LU_{ab}.\tag3.1
$$
\endproclaim

Set
$$
\Omega=\{(a,b,x)\in\Bp_1\times\Bp_1\times L\Bp_1: a,b\notin x(S^1)\}.
$$
Let $\frak H$ denote the space of those holomorphic cocycles $\frak 
h=(\frak h_{ab})_{a,b\in\Bp_1}$ of the covering $\{LU_a\}$, for which 
$\frak h_{ab}(x)$ depends holomorphically on $a,b$, and $x\in LU_{ab}$, and 
each
$\frak h_{ab}$ is cuspidal.
Then $\frak H\subset\Cal O(\Omega)$, with the compact open topology, is 
a complete, separated, locally convex space.
The action of $G$ on $\Omega$ induces a $G$--module structure on $\frak 
H$:
$$
(g^*\frak h)_{ab}(x)=\frak h_{ga,gb}(gx),\qquad g\in G.\tag3.2
$$
Proposition 3.2 implies the map $[f]\mapsto\frak f$ is a monomorphism 
$H^{0,1}(L\Bp_1)\to\frak H$ of $G$--modules.

\proclaim{Theorem 3.3}The map $[f]\mapsto \frak f$ is an isomorphism 
$H^{0,1}(L\Bp_1)\to\frak H$.
\endproclaim

The proof would be routine if the loop space $L\Bp_1$ admitted smooth
partitions of unity; but a typical loop space does not, see [K].
The proof that we offer here will work only when the loops in $L\Bp_1$ 
are of regularity $W^{1,3}$ at least, and we shall return to the case 
of $L_{1,p}(\Bp_1)$, $p<3$, in Section 6.

Those $g\in G$ that preserve the Fubini--Study metric form a subgroup 
(isomorphic to) SO$(3)$.
Denote the Haar probability measure on SO$(3)$ by $dg$.

\proclaim{Lemma 3.4}Unless $L\Bp_1=L_{1,p}\Bp_1,\ p<3$, there is a
$\chi\in C^\infty(L\Bp_1)$ such that $\chi=0$ in a neighborhood of
$L\Bbb P_1\setminus L\Bbb C=\{x\colon\infty\in x(S^1)\}$, and
$\int_{\SO(3)}g^*\chi dg= 1$.
\endproclaim

\demo{Proof}With $c_0\in (0,\infty)$ to be specified later, fix a 
nonnegative $\rho\in C^\infty(\Br)$ such that $\rho(\tau)=1$, resp.~0 when 
$|\tau|<c_0$, resp.~$>2c_0$.
For $x\in L\Bc$ let
$$
\psi(x)=\rho\left(\int_{S^1}(1+|x|^2)^{3/4}\right),
$$
and define $\psi(x)=0$ if $x\in L\Bp_1\backslash L\Bc$.
We claim that $\psi$ vanishes in a neighborhood of an arbitrary
$x\in L\Bp_1\setminus L\Bc$. This will then also imply that
$\psi\in C^\infty(L\Bp_1)$.

Indeed, suppose $x(t_0)=\infty$.
In a neighborhood of $t_0\in S^1$ the function $z=1/x$ is $W^{1,3}$, 
hence H\"older continuous with exponent 2/3 by the Sobolev Embedding 
Theorem, [H, Theorem 4.5.12].
In this neighborhood therefore $|x(t)|\geq c|t-t_0|^{-2/3}$, and
$\int_{S^1}(1+|x|^2)^{3/4}=\infty$.
When $y\in L\Bc$ is close to $x$, $\int_{S^1}(1+|y|^2)^{3/4}>2 c_0$, 
i.e.~$\psi(y)=0$.

Next we show that for every $x\in L\Bp_1$ there is a $g\in \SO(3)$ with 
$\psi(gx)>0$.
Let $d(a,b)$ denote the Fubini--Study distance between $a,b\in\Bp_1$; 
then with some $c>0$
$$
1+|\zeta|^2\leq {c\over d(\zeta,\infty)^2},\quad\text{ and }\quad
\int_{S^1} (1+|x|^2)^{3/4}\leq c\int_{S^1} d(x,\infty)^{-3/2}.
$$
Hence
$$
\int_{\SO(3)}\int_{S^1}(1+|gx(t)|^2)^{3/4}dt dg\leq 
c\int_{S^1}\int_{\SO(3)} d(gx(t),\infty)^{-3/2}dg dt=cI,
$$
where, for any $\zeta\in\Bp_1$
$$
I=\int_{\SO(3)} d(g\zeta,\infty)^{-3/2}dg=\int_{\Bp_1} 
d(\cdot,\infty)^{-3/2}<\infty,
$$
the last integral with respect to the Fubini--Study area form.
If $c_0$ is chosen $>cI$ then indeed $\int_{S^1}(1+|gx|^2)^{3/4}<c_0$ 
and
$\psi(gx)=1$ for some $g\in \SO(3)$.

It follows that $\int_{\SO(3)}\psi(gx)dg>0$, and we can take 
$\chi(x)=\psi(x)/\int_{\SO(3)}\psi(gx)dg$.
\enddemo

\demo{Proof of Theorem 3.3}Given $\frak h\in\frak H$, extend the 
functions $(g^*\chi)\frak h_{a,g\infty}$ from $LU_{a,g\infty}$ to $LU_a$ by 
zero, and define the cuspidal functions
$$
u_a=\int_{\SO(3)}(g^*\chi)\frak h_{a,g\infty}dg,\quad a\in\Bp_1.
$$
Then $u_a-u_b=\int_{\SO(3)}(g^*\chi)\frak h_{ab}dg=\frak h_{ab}$, so 
that $f=\overline\partial u_a$ on $LU_a$ consistently defines a closed 
$f\in C^\infty_{0,1}(L\Bp_1)$.
It is immediate that the map $\frak h\mapsto [f]\in H^{0,1}(L\Bp_1)$ is 
left inverse to the monomorphism $[f]\mapsto\frak f$, whence the 
theorem follows.
\enddemo

\heading{4.\ The map $\frak H\to\frak F$}\endheading

Consider an $\frak h=(\frak h_{ab})\in\frak H$.
The cocycle relation implies that $d_\zeta\frak h_{a\zeta}(x)$ is 
independent of $a$; for $\zeta\in\Bc$ we can write it as
$$
d_\zeta\frak 
h_{a\zeta}(x)=F\biggl(\zeta,{1\over\zeta-x}\biggr)d\zeta,\qquad x\in LU_\zeta,\tag4.1
$$
where $F\in\Cal O(\Bc\times L\Bc)$.
Set $F=\alpha(\frak h)$.
Since $\frak h_{aa}=0$,
$$
\frak h_{ab}(x)=\int_a^b 
F\biggl(\zeta,{1\over\zeta-x}\biggr)d\zeta,\tag4.2
$$
provided $a,b$ are in the same component of $\Bp_1\backslash 
x(S^1)$---that we shall express by saying $x$ does not separate $a,b$---, and we 
integrate along a path within this component.
The main result of this section is

\proclaim{Theorem 4.1} $\alpha(\frak h)=F\in\frak F$.
\endproclaim

The heart of the matter will be the special case when $\frak h$ is in 
an
irreducible submodule $\approx \frak K_n$.
A vector that corresponds in this isomorphism to 
const$(d\zeta)^{-n}\in\frak K_n$ is said to be of lowest weight $-n$.
Thus, if $\frak l$ is of lowest weight $-n\leq 0$, then
$$
\alignat 2
g^*_\lambda\frak l&=\lambda^{-n}\frak l,\quad
&&\text{when }g_\lambda\zeta=\lambda\zeta,\quad \lambda\in\Bc,\text{ 
and}\tag4.3\\
g^*_\lambda\frak l&=\frak l,\quad
&&\text{when }g_\lambda\zeta=\zeta+\lambda,\quad \lambda\in\Bc.\tag4.4
\endalignat
$$
Conversely, an $\frak l\not= 0$ satisfying (4.3), (4.4) is a lowest 
weight vector and spans an irreducible submodule, isomorphic to $\frak 
K_n$,
but we shall not need this fact.

If $\frak l\in\frak H$ satisfies (4.4) then $\frak 
l_{\infty\zeta}(x)=\frak l_{\infty,\zeta+\lambda}(x+\lambda)$ by (3.2), whence 
$d_\zeta\frak l_{\infty\zeta}(x)$ depends only on $\zeta-x$, and $\alpha(\frak l)$ 
is of form $F(\zeta,y)=E(y)$.
If, in addition, $\frak l$ satisfies (4.3), then similarly it follows 
that $E\in\Cal O(L\Bc)$ is homogeneous of degree $n+1$.
We now fix a nonzero lowest weight vector $\frak l\in\frak H$, the 
corresponding $(n+1)$--homogeneous polynomial $E$, and its polarization 
$\Cal E$, cf.~(1.2).

\proclaim{Proposition 4.2}$\Cal E(1\otimes y_1\otimes\ldots\otimes 
y_n)=0$, and so $E(y+\text{const})=E(y)$.
\endproclaim

\demo{Proof}Since $\frak l_{\infty 0}\in\Cal O(LU_{\infty 0})$ is 
cuspidal and homogeneous of order $-n$,
$$
0=\lim_{\lambda\to\infty}\frak l_{\infty 0}
\biggl({1\over\lambda+x}\biggr)=\lim_{\lambda\to\infty}
\lambda^n\frak l_{\infty 0}\biggl({1\over 1+x/\lambda}\biggr).
$$
Thus $\frak l_{\infty 0}$ vanishes at 1 to order $\geq n+1$.
Hence
$$
{\partial\over\partial\zeta}\bigg|_{\zeta=0} \frak l_{\infty 
0}(x-\zeta)={\partial\over\partial\zeta}\bigg|_{\zeta=0}\frak 
l_{\infty\zeta}(x)=E\biggl({1\over x}\biggr)
$$
vanishes at $x=1$ to order $\geq n$, and the same holds for $E(x)$.
Differentiating $E$ in the directions $y_1,\ldots,y_n$, we obtain at 
$x=1$, as needed, that $n!\Cal E(1\otimes y_1\otimes\ldots\otimes 
y_n)=0$.
\enddemo

Let $\frak K_n\ni\varphi\mapsto\frak h^\varphi\in\frak H$ denote the 
homomorphism that maps $(d\zeta)^{-n}$ to $\frak l$.

\proclaim{Proposition 4.3}
$$
d_\zeta\frak h_{a\zeta}^\varphi(x)=\psi(\zeta)
E\biggl({1\over\zeta-x }\biggr)d\zeta,
\qquad\varphi(\zeta)=\psi(\zeta)(d\zeta)^{-n}.\tag4.5
$$
By homogeneity, the right hand side can also be written 
$\varphi(\zeta)E(d\zeta/(\zeta-x))$.
\endproclaim

\demo{Proof}Denote the form on the left hand side of (4.5) by 
$\omega^\varphi$.
In view of (3.2) it transforms under the action of $G$ on
$\Bp_1\times L\Bp_1$ as
$$
g^*\omega^\varphi=\omega^{g\varphi},\qquad g\in G.\tag4.6
$$
If we show that the right hand side of (4.5) transforms the same way, 
then (4.5) will follow, since it holds when $\psi\equiv 1$, see (4.1).
In fact, it will suffice to check the transformation formula for 
$g\zeta=\lambda\zeta,\ g\zeta=\zeta+\lambda\ (\lambda\in\Bc)$, and 
$g\zeta=1/\zeta$, maps that generate $G$.
We shall do this for the last map, the most challenging of the three 
types.
The pullback of the right hand side of (4.5) by $g\zeta=1/\zeta$ is
$$
\align
(g\varphi)(\zeta)E\bigg({d(g\zeta)\over 
g\zeta-gx}\bigg)&=(g\varphi)(\zeta)E\bigg({-d\zeta/\zeta^2\over
(1/\zeta)-(1/x)}\bigg)\\
&=(g\varphi)(\zeta)E\bigg({d\zeta\over\zeta-x}-
{d\zeta\over\zeta}\bigg)=(g\varphi)(\zeta)E\bigg({d\zeta\over\zeta-x}\bigg),
\endalign
$$ by Proposition 4.2, which is what we needed.
\enddemo

The form $\Cal E$ defines a symmetric distribution $D$ on the torus 
$T=(S^1)^{n+1}$ as in Section 1, cf.~(1.14).
By (1.15), (4.2), and Proposition 4.3
$$
\frak h^\varphi_{ab}(x)=\int_a^b\psi(\zeta)\biggl\langle D,\ 
{1\over\zeta-x}\otimes\ldots\otimes{1\over\zeta-x}\biggr\rangle d\zeta,
\qquad \varphi=\psi(\zeta)(d\zeta)^{-n},\tag4.7
$$
provided $x\in L_\infty U_{ab}$ does not separate $a,b$.
To prove Theorem 4.1, we have to understand supp $D$.
Let
$$
O=\{x\in C^\infty(S^1)\colon\pm i\not\in x(S^1)\},\text{ and }O'=\{x\in 
O\colon[-i,i]\cap x(S^1)=\emptyset\},
$$
where $[-i,i]$ stands for the segment joining $\pm i$.

\proclaim{Lemma 4.4}With $\Delta$ a symmetric distribution on
$T=(S^1)^{n+1}$ and $\nu=0,\ldots,2n-2$, let
$$
I_\nu(x)=\int_{[-i,i]}\biggl\langle\Delta,\ 
{1\over\zeta-x}\otimes\ldots\otimes{1\over\zeta-x}\biggr\rangle\zeta^\nu d\zeta,\qquad x\in O'.
$$
If each $I_\nu$ continues analytically to $O$ then $\Delta$ is 
supported on the diagonal of $T$.
\endproclaim

In preparation to the proof, consider a holomorphic vector field $V$ on 
$O$, and observe that $VI_\nu$ also continues analytically to $O$.
Such vector fields can be thought of as holomorphic maps $V\colon O\to 
C^\infty(S^1)$.
Using the symmetry of $\Delta$ we compute
$$
(VI_\nu)(x)=(n+1)\int_{[-i,i]}\biggl\langle\Delta,\ {V(x)\over 
(\zeta-x)^2}\otimes 
{1\over\zeta-x}\otimes\ldots\otimes{1\over\zeta-x}\biggr\rangle \zeta^\nu d\zeta,\quad x\in O'.\tag4.8
$$

\def\os{\overline s}
\def\var{\varepsilon}
\demo{Proof of Lemma 4.4, case $n=1$}Let $\overline s_0\not= \overline 
s_1\in S^1$.
To show $\Delta$ vanishes near $\os=(\os_0,\os_1)$, construct a smooth 
family
$x_{\var,s}\in O$ of loops, where $\var\in [0,1]$ and $s\in T$ is in a
neighborhood of $\os$, so that
$$
x_{\var,s}(\tau)=(-1)^j(\var^2+(\tau-s_j)^2),
\qquad\text{when }\tau\in S^1\text{ is near }\os_j,\ j=0,1;\tag4.9
$$
here, perhaps abusively, $\tau-s_j$ denotes both a point in
$S^1=\Bbb R/\Bbb Z$ and its representative in $\Bbb R$ that is closest 
to
$0$. Make sure that $x_{\var,s}\in O'$ when $\var > 0$.
Fix $y_0,y_1\in C^\infty(S^1)$ so that $y_j\equiv 1$ near $\os_j$, and 
(4.9) holds when $\tau, s_j$ are in a neighborhood of supp $y_j$.
This forces $y_0,y_1$ to have disjoint support.
With constant vector fields $V_j=y_j$
$$
(V_1 V_0 I_0)(x)=2\int_{[-i,i]}\biggl\langle\Delta,\ {y_0\over 
(\zeta-x)^2}\otimes{y_1\over (\zeta-x)^2}\biggr\rangle d\zeta,\qquad x\in 
O',\tag4.10
$$
analytically continues to $O$.
In particular, for $\var>0$ and $t=(t_0,t_1)\in T$ setting
$$
K_\var(t,s)=\int_{[-i,i]}\ {y_0(t_0)y_1(t_1)d\zeta\over 
(\zeta-x_{\var,s}(t_0))^2(\zeta-x_{\var,s}(t_1))^2},\qquad s\text{ near }\os,
$$
it follows that $\langle\Delta,K_\var(\cdot,s)\rangle$ stays bounded as 
$\var\to 0$.
Therefore, if $\rho\in C^\infty(T)$ is supported in a sufficiently 
small neighborhood of $\os$,
$$
\langle\Delta,\ \var^4\int_T K_\var (\cdot,s)\rho(s)ds\rangle\to 
0,\qquad \var\to 0.\tag4.11
$$
On the other hand we shall show that for such $\rho$
$$
\var^4\int_T K_\var (\cdot,s)\rho(s)ds\to c\rho,\qquad\var\to 
0,\tag4.12
$$
in the topology of $C^\infty(T)$; here $c\not= 0$ is a constant.

It will suffice to verify (4.12) on supp $y_0\otimes y_1$, since both 
sides
vanish on the complement.
Thus we shall work on small neighborhoods of $\os$; we can pretend 
$\os\in\Br^2$, and work on $\Br^2$ instead of $T$.
When $s,\ t\in\Br^2$ are close to $\os$, the left hand side of (4.12) 
becomes
$$
\var^4 y_0 (t_0)y_1(t_1)\int_{\Br^2}\int_{[-i,i]}\
{\rho(s)d\zeta ds\over (\zeta-\var^2-(s_0-t_0)^2)^2
(\zeta+\var^2+(s_1-t_1)^2)^2}.\tag4.13
$$
Substituting $s=t+\var u$ and $\zeta=\var^2\xi$, we compute the limit 
in (4.12) is
$$
\gathered
\lim_{\var\to 0} 
y_0(t_0)y_1(t_1)\int_{\Br^2}\int_{[-i/\var^2,i/\var^2]}\
{\rho(t+\var u)d\xi du\over (\xi-1-u_0^2)^2(\xi+1+u_1^2)^2 }\\
=4\pi i y_0(t_0) y_1 (t_1)\ \int_{\Br^2}\ {\rho(t) du\over
(2+u_0^2+u_1^2)^3}=c\rho(t),
\endgathered\tag4.14
$$
if $y_0\otimes y_1=1$ on supp$\,\rho$.
This limit is first seen to hold uniformly.
However, since the integral operator in (4.13) is a convolution, in 
(4.14) in fact all derivatives converge uniformly.
Now (4.11) and (4.12) imply $\langle\Delta,\rho\rangle=0$, so that 
$\Delta$ vanishes close to $\os$, q.e.d.
\enddemo

\demo{Proof of Lemma 4.4, general $n$}The base case $n=1$ settled and 
the statement being vacuous when $n=0$, we prove by induction.
Assume the Lemma holds on the $n$--dimensional torus, and with $y\in 
C^\infty(S^1)$, consider holomorphic vector fields $V_\mu(x)=yx^\mu,\ 
\mu=0,1,2$.
(These vector fields continue to all of $L\Bp_1$, and generate the Lie 
algebra of the loop group $LG$.)
In view of (4.8), for $x\in O'$
$$
\int_{[-i,i]}\biggl\langle\Delta,y\otimes 
{1\over\zeta-x}\otimes\ldots\otimes{1\over\zeta-x}\biggr\rangle\zeta^\nu d\zeta={1\over n+1}(V_0 
I_{\nu+2}-2 V_1 I_{\nu+1}+V_2 I_\nu).\tag4.15
$$
Therefore the left hand side continues analytically to $O$, provided 
$\nu=0,\ldots,2n-4$.
If $\Delta^y$ denotes the distribution on $(S^1)^n$ defined by 
$\langle\Delta^y,\rho\rangle=\langle\Delta,y\otimes\rho\rangle$, the left hand 
side of (4.15) is
$$
\int_{[-i,i] }\biggl\langle\Delta^y,{1\over\zeta-x}\otimes\ldots\otimes 
{1\over\zeta-x}\biggr\rangle\zeta^\nu d\zeta.
$$
The inductive hypothesis implies $\Delta^y$ is supported on the 
diagonal of
$(S^1)^n$.
This being true for all $y$, the symmetric distribution $\Delta$ itself
must be supported on the diagonal.
\enddemo

\proclaim{Corollary 4.5}The distribution $D$ in (4.7) is supported on 
the diagonal of $T$.
\endproclaim

\demo{Proof of Theorem 4.1}First assume that $\frak h\in\frak H$ is in 
an irreducible submodule $\approx \frak K_n$, and $\frak l\not= 0$ is a 
lowest weight vector in this submodule.
Thus $\frak h=\frak h^\varphi$ with some $\varphi\in\frak K_n$, 
$\varphi(\zeta)=\psi(\zeta)(d\zeta)^{-n}$.
With $\frak l$ we associated an $(n+1)$--homogeneous polynomial $E$ on 
$L\Bc$ and a distribution $D$ on $(S^1)^{n+1}$.
By Proposition 4.3 $F(\zeta,y)=\psi(\zeta)E(y)$, and so 
$F(\zeta,y+\text{const})=F(\zeta,y)$ by Proposition 4.2.
Since deg $\psi\leq 2n,\ F(\zeta/\lambda,\lambda^2 y)=O(\lambda^2)$ as 
$\lambda\to 0$.
Finally, take $x,y\in L\Bc$ with disjoint support.
If $x,y\in C^\infty(S^1)$,
$$
E(x+y)=\langle D,(x+y)^{\otimes n+1}\rangle=\langle D, x^{\otimes 
n+1}\rangle+\langle D,y^{\otimes n+1}\rangle=E(x)+E(y),
$$
as supp $D$ is on the diagonal.
By approximation $E(x+y)=E(x)+E(y)$ follows in general, whence $F$ 
itself is additive.
We conclude $F\in\frak F$ if $\frak h$ is in an irreducible submodule.

By linearity it follows that $F\in\frak F$ whenever $\frak h$ is in the 
span of irreducible submodules.
Since this span is dense in $\frak H$ (cf.~[BD, III.5.7] and the 
explanation in the Introduction connecting representations of $G$ with those 
of the compact group SO$(3)$), $\alpha(\frak h)\in\frak F$ for all 
$\frak h\in\frak H$.
\enddemo

\proclaim{Theorem 4.6}The map $\alpha$ is a $G$--morphism.
\endproclaim

\demo{Proof}It suffices to verify that the restriction of $\alpha$ to 
an irreducible submodule of $\frak H$ is a $G$--morphism, which follows 
directly from Proposition 4.3.
\enddemo

\heading 5. The structure of $\frak H$\endheading

The main result of this Section is

\proclaim{Theorem 5.1}The $G$--morphism $\alpha\colon\frak H\to\frak F$ 
has a right inverse $\beta$.
Its kernel is one dimensional, spanned by the $G$--invariant cocycle
$$
\frak h_{ab}(x)=ind_{ab} x
$$
(= the winding number of $x\colon S^1\to U_{ab})$.
\endproclaim

We shall need the following

\proclaim{Lemma 5.2}With notation as in Section 1, suppose 
$z_1,\ldots,z_N\in L^-\Bc$ are such that no point in $S^1$ is contained in the 
support of more than two $z_j$.
If $\tilde F\in\tilde\frak F$ then
$$
\tilde F(\zeta,\sum^N_{j=1} z_j)=\sum_{i<j}\tilde 
F(\zeta,z_i+z_j)-(N-2)\sum^N_{j=1}\tilde F(\zeta,z_j).\tag5.2
$$
In particular, if $N\ge 3$, and writing $z_0=z_N$ only consecutive
$\supp\, z_j$'s intersect each other, then
$$
\tilde F(\zeta,\sum_1^N z_j)=
\sum_1^N \tilde F(\zeta, z_{j-1}+z_j)-\sum_1^N\tilde F(\zeta,z_j).
$$
\endproclaim

\demo{Proof}It will suffice to verify (5.2) when
$\tilde F(\zeta,z)=\tilde E(z)$ is homogeneous, in which case it 
follows by expressing both sides in terms of the polarization of $\tilde E$, 
and
using Lemma 1.2a. The second formula follows from (5.2) by
applying additivity to terms with non consecutive $i,j$.
\enddemo

\demo{Proof of Theorem 5.1}(a) Construction of the right inverse.
By Theorem 1.1, for $F\in\frak F$ we can choose $\tilde F\in\tilde\frak 
F$, depending linearly on $F$, so that $F(\zeta,y)=\tilde F(\zeta,\dot 
y)$.
With $x\in L\Bp_1$ consider the differential form
$$
F\biggl(\zeta,{1\over\zeta-x}\biggr)d\zeta=\tilde F\biggl(\zeta,
{\dot x\over (\zeta-x)^2}\biggr)d\zeta,\tag5.3
$$
holomorphic in $\Bc\backslash x(S^1)$.
In fact, it is holomorphic at $\zeta=\infty$ as well, provided 
$\infty\not\in x(S^1)$, since the coefficient of $d\zeta$ vanishes to second 
order at $\zeta=\infty$.
This latter is easily verified when $\tilde F(\zeta,z)=\zeta^\nu\tilde 
E(z)$ and $\tilde E$ is $(n+1)$--homogeneous, $\nu\leq 2n$; in general 
it follows from the posthomogeneous expansion
$$
\tilde F(\zeta,z)=\sum^\infty_0\tilde 
F_n(\zeta,z)=\sum^\infty_{n=0}\sum^{2n}_{\nu=0}\zeta^\nu\tilde E_{n\nu}(\zeta).
$$

Hence, if $x\in L\Bp_1$ does not separate $a$ and $b$, the integral
$$
h_{ab}(x)=\int_a^b\tilde F\biggl(\zeta,\ {\dot x\over 
(\zeta-x)^2}\biggr)d\zeta\tag5.4
$$
is independent of the path joining $a,b$ within $\Bp_1\backslash 
x(S^1)$, and defines a holomorphic function of $a,b$, and $x$.

We claim that $h_{ab}$ can be continued to a cuspidal cocycle
$\frak h=(\frak h_{ab})\in\frak H$.
First we prove a variant.
Let $\sigma\in C^\infty(S^1)$ be supported in a closed arc $I\not= 
S^1$.
Given finitely many $a,b,\ldots\in\Bp_1$, set
$$
W_{ab\ldots}=\{x\in L\Bp_1\colon a,b,\ldots\not\in x (I)\}\supset 
LU_{ab\ldots}.
$$
We shall show that the integrals
$$
\int_a^b\tilde F\biggl(\zeta,\ {\sigma\dot x\over 
(\zeta-x)^2}\biggr)d\zeta,
\qquad x\text{ does not separate }a,b,\tag5.5
$$
can be continued to functions $\frak k_{ab}(x)$ depending 
holomorphically on $a,b\in\Bp_1$, and $x\in W_{ab}$.
The main point will be that, unlike $LU_{ab\ldots}$, the sets 
$W_{ab\ldots}$ are connected.

If $x_1\in W_{ab}$, construct a continuous curve $[0,1]\ni\tau\mapsto 
x_\tau\in W_{ab}$, $x_0=$ constant loop.
Cover $S^1$ with open arcs $J_1,\ldots,J_N=J_0$, $N\ge 3$,
so that only
consecutive $\overline J_j$'s intersect, and no $x_\tau(\overline 
J_i\cup\overline J_j)$ separates $a$ and $b$.
Choose a $C^\infty$ partition of unity $\{\rho_j\}$ subordinate to 
$\{J_j\}$.
For $x$ in a connected neighborhood $W\subset W_{ab}$ of 
$\{x_\tau\colon 0\leq\tau\leq 1\}$ define
$$
\frak k_{ab}(x)=
\sum_1^N\int_a^b\tilde F\biggl(\zeta,\ {(\rho_{j-1}+\rho_j)\sigma\dot x
\over (\zeta-x)^2}\biggr)d\zeta-\sum^N_1\int_a^b\tilde F
\biggl(\zeta,\ {\rho_j\sigma\dot x\over 
(\zeta-x)^2}\biggr)d\zeta.\tag5.6
$$
In the first sum we extend
$(\rho_{j-1}+\rho_j)\sigma\dot x/(\zeta-x)^2$ to
$S^1\backslash (J_{j-1}\cup J_j)$ by $0$, and integrate along paths in 
$\Bp_1\backslash x(\overline J_{j-1}\cup \overline J_j)$; we interpret
the second sum similarly.
The neighborhood $W$ is to be chosen so small that no $x(\overline 
J_i\cup\overline J_j)$ separates $a$ and $b$ when $x\in W$.

As above, the integrals in (5.6) are independent of the path, and 
define a holomorphic function in $W$.
By Lemma 5.2, $\frak k_{ab}$ agrees with (5.5) when $x$ is near $x_0$.
Furthermore, the germ of $\frak k_{ab}$ at $x_1$ depends on the curve 
$x_\tau$ only through the choice of the $\rho_j$.
In fact, it does not even depend on $\rho_j$:\ let $\frak k'_{ab}$ be 
the function obtained if in (5.6) the $\rho_j$ are replaced by another 
partition of unity $\rho'_h$.
It will suffice to show that $\frak k_{ab}=\frak k'_{ab}$ under the 
additional assumption that each $\rho'_h$ is supported in some $J_j$.
In this case $\frak k'_{ab}$ is holomorphic in $W$ and agrees with 
$\frak k_{ab}$ near $x_0$, hence on all of $W$.

Therefore, by varying the partition of unity $\rho_j$, we can use (5.6) 
to define $\frak k_{ab}(x)$ depending holomorphically on $a,b\in\Bp_1$, 
and $x\in W_{ab}$.
Also, $\frak k_{ab}+\frak k_{bc}=\frak k_{ac}$ on $W_{abc}$, since this 
is so in a neighborhood of constant loops, and $W_{abc}$ is connected.

Now, to obtain a continuation of $h_{ab}$ in (5.4), construct a 
partition of unity $\sigma_1,\sigma_2,\sigma_3\in C^\infty(S^1)$, so that
$\supp\,(\sigma_i+\sigma_j)\neq S^1$ and 
$\bigcap_1^3\supp\,\sigma_j=\emptyset$.
Setting $\sigma_0=\sigma_3$, in light of Lemma 5.2 we can rewrite (5.4)
$$
h_{ab}(x)=\sum_1^3\int_a^b\tilde F\biggl(\zeta,
{(\sigma_{j-1}+\sigma_j)\dot x
\over (\zeta-x)^2 }\biggr)d\zeta-
\sum_1^3\int_a^b\tilde F\biggl(\zeta,{\sigma_j\dot x\over(\zeta-x)^2 
}\biggr)d\zeta,
$$
and continue each term to $LU_{ab}$, as above.
We obtain a holomorphic cocycle $\beta(F)=\frak h=(\frak h_{ab})$, with 
$\frak h_{ab}$ depending holomorphically on $a,b$, and one easily 
checks that each $\frak h_{ab}$ is cuspidal.
Therefore $\beta(F)\in\frak H$.
Finally, $\alpha\beta(F)$ can be computed by considering $d_\zeta\frak 
h_{a\zeta}(x)$ with $a$ in the same component of $\Bp_1\backslash 
x(S^1)$ as $\zeta$, so that (5.4) gives
$$
d_\zeta\frak h_{a\zeta}(x)=d_\zeta h_{a\zeta}(x)=\tilde 
F\biggl(\zeta,{\dot x\over (\zeta-x)^2}\biggr)d\zeta=F\biggl(\zeta,{1\over 
\zeta-x}\biggr)d\zeta.
$$
Thus $\alpha\beta(F)=F$ as needed.

(b) The kernel of $\alpha$.
Take an irreducible submodule of Ker $\alpha$, spanned by a vector
$\frak l$ of lowest weight $-n\leq 0$.
Since $F=\alpha(\frak l)=0$, (4.2) implies $\frak l_{ab}(x)=0$ if $x$ 
does not separate $a,b$; hence, by analytic continuation, whenever 
ind$_{ab}x=0$.
By the cocycle relation $\frak l_{ac}(x)=\frak l_{bc}(x)$ if 
ind$_{ab}x=0$, i.e., if ind$_{ac}x=\text{ ind}_{bc}x$.

Consider the components of $LU_{0\infty}$
$$
X_r=\{x\in LU_{0\infty}\colon\text{ ind}_{0\infty} x=r\},\qquad 
r\in\Bz.
$$
Let
$$
x_1(t)=e^{irt},\quad y(t)=e^{2irt}+e^{3irt-4}.\tag5.7
$$
We shall shortly show that whenever $x\in LU_{0\infty}$ is in a 
sufficiently small neighborhood of $x_1$, and 
$(\kappa,\lambda)\in\Bc^2\backslash (0,0)$, then $z_{\kappa\lambda}=\kappa x+\lambda y\in X_r+\Bc$.
It follows that with such $x,y$ we can define $h(\kappa,\lambda)=\frak 
l_{a\infty}(z_{\kappa\lambda})$, where $a$ is chosen so that 
$\text{ind}_{a\infty}z_{\kappa\lambda}=r$.
Thus $h\in\Cal O(\Bc^2\backslash (0,0))$, and by Hartogs' theorem it
extends to all of $\Bc^2$; also, it is homogeneous of degree $-n$.
It follows that $h$ is constant, indeed zero when $n>0$.
In all cases $\frak l_{0\infty}(x)=h(1,0)=h(0,1)$ is independent of 
$x$.
This being true for $x$ in a nonempty open set, $\frak l_{0\infty}$ is 
constant on $X_r$.
It follows that $\frak l_{a\infty}(x)=\frak l_{0\infty}(x-a)$ is 
locally constant, and so is $\frak l_{ab}=\frak l_{a\infty}-\frak 
l_{b\infty}$.
Moreover, $\frak l_{ab}=0$ unless $n=0$.

Suppose now $n=0$, and let $\frak l_{0\infty}|X_1=l\in\Bc$.
We have $\frak l_{a\infty}(x)=\frak l_{0\infty}(x-a)=l$ if 
ind$_{a\infty}x=1$.
Choose a homeomorphic $x\in L\Bc$ and $a,b\in\Bc\backslash x(S^1)$ so 
that ind$_{ab}x=1$; say $b$ is in the unbounded component.
Then $\frak l_{ab}(x)=\frak l_{a\infty}(x)-\frak l_{b\infty}(x)=l$, and 
the same will hold if $x$ is slightly perturbed.
It follows that $\frak l_{ab}(x)=l$ whenever ind$_{ab}x=1$, and in this 
case $\frak l_{ba}(x)=-l$.
Finally, with a generic $y\in LU_{ab}$ choose $a_0=a,a_1,\ldots,a_m=b$ 
in $\Bp_1\backslash y(S^1)$ so that ind$_{a_{j-1}a_j}y=\pm 1$.
Then
$$
\frak l_{ab}(y)=\sum_1^m\frak l_{a_{j-1}a_j }(y)=l\sum_1^m\text{ 
ind}_{a_{j-1}a_j}y=l\text{ ind}_{ab}y.
$$
The upshot is that any irreducible submodule of Ker $\alpha$ is spanned 
by $\frak h$ in (5.1), whence Ker $\alpha$ itself is spanned by $\frak 
h$, as claimed.

We still owe the proof that $\kappa x+\lambda y\in X_r+\Bc$ unless 
$\kappa=\lambda=0$, for $x$ near $x_1$ and $y$ given in (5.7).
In fact, the general statement follows once we prove it for $r=1$ and 
$x=x_1$, that we henceforward assume.
If $|\kappa|\geq 2|\lambda|$ then $z_{\kappa\lambda}\in X_1$ by 
Rouch\'e's theorem.
Otherwise consider the polynomial
$$
P(\zeta)=\kappa\zeta+\lambda(\zeta^2+e^{-4}\zeta^3),\qquad\zeta\in\Bc.
$$
For fixed $|\zeta|<2$ the equation $P(\eta)=P(\zeta)$ has two solutions 
with $|\eta|<5$, again by Rouch\'e's theorem.
One of the solutions is $\eta=\zeta$.
Let $\eta=R(\zeta)$ be the other one, so that $R$ is holomorphic.
There are only finitely many $\zeta$ with $|\zeta|=|R(\zeta)|=1$.
Indeed, otherwise $|R(\zeta)|=1$ would hold for all unimodular $\zeta$, 
and by the reflection principle $R$ would be rational.
However, $P(R(\zeta))=P(\zeta)$ cannot hold with rational 
$R(\zeta)\not=\zeta$.
We conclude that $z_{\kappa\lambda}(S^1)$ has only finitely many 
self--intersection points.

Since $P(0)=0,\text{ ind}_{0\infty}z_{\kappa\lambda}\geq 1$.
Drag a point $a$ from 0 to $\infty$ along a path that avoids multiple 
points of $z_{\kappa\lambda}(S^1)$.
Each time we cross $z_{\kappa\lambda}(S^1)$, 
ind$_{a\infty}z_{\kappa\lambda}$ changes by $\pm 1$.
It follows that ind$_{a\infty}z_{\kappa\lambda}=1$ for some $a$, which 
completes the proof.
\enddemo

For the space $L_{1,p}\Bp_1$ Theorems 2.1, 2.2, and the construction in 
Theorem 5.1 lead to explicit representations of elements of $\frak H$.
First there are the multiples of the cocycle (5.1), and then there is 
the complementary subspace $\beta(\frak F)=\oplus_{n\leq p-1}\beta(\frak 
F^n)$, see Theorem 2.2.
According to Theorems 2.1, 2.2 elements of $\frak F^n$ are of form
$$
F(\zeta,y)=\sum^{2n}_{\nu=0}\zeta^\nu\langle\Phi_\nu,\dot 
y^{n+1}\rangle,\qquad \Phi_\nu\in L^{p/(n+1)}(S^1)^*.
$$
Following the proof of Theorem 5.1, to compute $\frak h=\beta(F)$ we
set $\tilde 
F(\zeta,z)=\sum_\nu\zeta^\nu\langle\Phi_{\nu},z^{n+1}\rangle$.
The substitution $\zeta=\xi+c$ shows that
$$
R_\nu(a,b,c)=\int_a^b {\zeta^\nu d\zeta\over (\zeta-c)^{2n+2}},\qquad 
0\leq\nu\leq 2n,\ c\in\Bp_1\backslash \{a,b\},
$$
are rational functions with poles at $c=a,b$, so that
$$
\frak h_{ab}(x)=\int_a^b\tilde F\biggl(\zeta,{\dot x\over 
(\zeta-x)^2}\biggr)d\zeta=\sum^{2n}_{\nu=0}\langle\Phi_\nu, R_\nu(a,b,x)\dot 
x^{n+1}\rangle,
$$
when $x$ does not separate $a,b$.
However, the right hand side makes sense for any $x\in LU_{ab}$ and, as 
one checks, defines $\frak h=\beta(F)$.
For example, if $F$, hence $\frak h$ are of lowest weight, then 
$\Phi_\nu=0$ for $\nu\geq 1$, and
$$
\frak h_{ab}(x)=\biggl\langle\Phi_0,\ {\dot x^{n+1}\over 2n+1}\
\biggl({1\over (x-a)^{2n+1}}-
{1\over (x-b)^{2n+1}}\biggr)\biggr\rangle.\tag5.8
$$
Letting $n=0$ and $\langle\Phi_0,z\rangle=\int_{S^1}z/2\pi i$, formula 
(5.8) recovers the locally constant cocycle (5.1) as well.
Thus we proved

\proclaim{Theorem 5.3}In the case of $W^{1,p}$ loop spaces, any lowest 
weight cocycle in the $n$'th isotypical subspace $\frak H^n\subset\frak 
H$ is of form (5.8) with (a unique) $\Phi_0\in L^{p/(n+1)}(S^1)$, 
$0\leq n\leq p-1$.
\endproclaim

\heading 6.\ Synthesis\endheading

In this last section we show how the results obtained by now imply the 
theorems of the Introduction.
Theorems 0.1 and 0.2 follow from the isomorphism
$H^{0,1}(L\Bp_1)\approx\frak H$ of $G$--modules (Theorem 3.3) and from 
the isomorphism $\frak H\approx\Bc\oplus\frak F$, a consequence of 
Theorem 5.1.
In particular, $H^{0,1}(L\Bp_1)^G\approx \Bc\oplus\frak F^0$.
The latter being isomorphic to the dual of $L^-\Bc=C^{k-1}(S^1)$, 
resp.~$W^{k-1,p}(S^1)$ by Theorem 2.1, Theorem 0.3 also follows.
Finally, Theorem 0.4 is a consequence of Theorems 2.2 and 2.1.

Seemingly we are done with all the proofs.
However, Theorem 3.3 has not yet been proved for loop spaces 
$L_{1,p}\Bp_1,\ p<3$, and we still have to revisit spaces of loops of low 
regularity.
This will give us the opportunity to explicitly represent classes in 
$H^{0,1}(L_{1,p}\Bp_1)$, in fact, for all $p\in [1,\infty)$.

Generally, given a complex manifold $M$, $1\leq p<\infty$, and a 
natural number $m\leq p$, consider the space $C^\infty_{0,q}((T^*M)^{\otimes 
m})$ of $(T^*M)^{\otimes m}$ valued $(0,q)$ forms on $M$.
If $\omega$ is such a form, $v\in\oplus^q T_s^{0,1}M$, and $w\in 
T_s^{1,0}M$, we can pair $\omega(v)\in (T_s^* M)^{\otimes m}$ with 
$w^{\otimes m}$, to obtain what we shall denote $\omega(v,w^m)\in\Bc$.
Write $LM$ for the space of $W^{1,p}$ loops in $M$, and observe that 
the tangent space $T_x^{0,1}LM$ is naturally isomorphic to the space
$W^{1,p}(x^*T^{0,1}M)$ of $W^{1,p}$ sections of the induced bundle $x^* 
T^{0,1}M\to S^1$ (see [L2, Proposition 2.2]  in the case of $C^k$ 
loops).

There is a bilinear map
$$
I=I_q\colon L^{p/m} (S^1)^*\times C^\infty_{0,q} ((T^* M)^{\otimes 
m})\to C^\infty_{0,q}(LM),
$$
obtained by the following Radon type transformation.
If
$$
(\Phi,\omega)\in L^{p/m}(S^1)^*\times C_{0,q}^\infty((T^*M)^{\otimes 
m}),
$$
$x\in LM$, and $\xi\in\oplus^q T_x^{0,1} LM\approx \oplus^q W^{1,p}(x^* 
T^{0,1}M)$, then $\omega(\xi,\dot x^m)\in L^{p/m}(S^1)$.
Define $I(\Phi,\omega)=f$ by
$$
f(\xi)=\langle\Phi,\omega(\xi,\dot x^m)\rangle.
$$
One verifies that $\overline\partial 
I(\Phi,\omega)=I(\Phi,\overline\partial\omega)$, whence $I_q$ induces a bilinear map
$$
L^{p/m}(S^1)^*\times H^{0,q}((T^*M)^{\otimes m})\to H^{0,q}(LM).
$$

Henceforward we take $M=\Bp_1,\ q=1,\ m=n+1$, and $\omega$ given on 
$\Bc$ by
$$
\omega={-1\over 2n+1}\ 
{\overline\zeta^{2n}d\overline\zeta\otimes(d\zeta)^{n+1}\over (1+|\zeta|^{4n+2})^{(2n+2)/(2n+1)}},\qquad \zeta\in\Bc,
$$
so that $f=I_1(\Phi,\omega)$ is a closed form on $L\Bp_1$. Explicitly,
$$
f(\xi)={-1\over 2n+1}\biggl\langle\Phi,\ {\xi\overline x^{2n}\dot 
x^{n+1}
\over (1+|x|^{2n+2})^{(2n+2)/(2n+1)}}\biggr\rangle,\qquad\xi
\in T_x^{0,1}L\Bp_1.\tag6.1
$$
To compute its image in $\frak H$ under the map of Theorem 3.3, let
$$
\theta_a={1\over 2n+1}\ \biggl({\zeta^{-2n-1}\over
(1+|\zeta|^{4n+2})^{1/(2n+1)}}-\zeta^{-2n-1}+
(\zeta-a)^{-2n-1}\biggl)(d\zeta)^{n+1}\qquad\text{ on }U_a.
$$
Thus $\overline\partial\theta_a=\omega|U_a$, and the cuspidal functions 
$u_a=I_0(\Phi,\theta_a)\in C^\infty (LU_a)$ solve $\overline\partial 
u_a=f|LU_a$.
Hence the image of $f$ in $\frak H$ is
$$
\frak h_{ab}(x)=u_a (x)-u_b(x)=\biggl\langle\Phi,\ {\dot x^{n+1}\over 
2n+1}\
\biggl({1\over (x-a)^{2n+1} }-
{1\over (x-b)^{2n+1} }\biggr)\biggr\rangle.
$$
Comparing this with Theorem 5.3 we see that by associating with a 
lowest weight $\frak h\in\frak H^n$ the functional $\Phi=\Phi_0$ of (5.8), 
and then $f\in C^\infty_{0,1}(L\Bp_1)$ of (6.1), the image of $f$ in 
$\frak H$ will be $\frak h$.
In particular, the class $[f]\in H^{0,1}(L\Bp_1)$ is also of lowest 
weight $-n$.
Therefore the linear map $\frak h\mapsto [f]$, defined for
$\frak h\in\frak H^n$ of lowest weight, can be extended to a 
$G$--morphism $\frak H^n\to H^{0,1}(L\Bp_1)$, and then to a $G$--morphism
$\bigoplus_{n\le p-1}\frak H^n=\frak H\to H^{0,1}(L\Bp_1)$, inverse to 
the morphism $H^{0,1}(L\Bp_1)\to \frak H$ of Theorem 3.3.
This completes the proof of Theorem 3.3, and now we are really done.

\enddocument